\documentclass[12pt,a4paper]{article}
\usepackage{amsmath,mathrsfs,amsmath,amssymb,enumerate}
\usepackage[colorlinks=true,linkcolor=black,anchorcolor=black,citecolor=black,baseurl=empty,urlcolor=black]{hyperref}
\usepackage{cleveref}
\usepackage{dsfont}
\usepackage{comment}
\usepackage{indentfirst}
\usepackage{url}
\usepackage[pagewise,mathlines]{lineno}
\numberwithin{equation}{section}
\newtheorem{thm}{Theorem}[section]

\newtheorem{prop}[thm]{Proposition}
\newtheorem{defn}[thm]{Definition}
\newtheorem{rem}[thm]{\bf Remark}

\mathchardef\lz="2D

\newcommand{\R }{\mathbb{R}^N }

\newcommand{\HO}{{H}_0^{1}(\varOmega )}

\newcommand{\RR}{\mathbb{R}}

\allowdisplaybreaks

\crefname{equation}{Problem}{Problems}
\crefname{thm}{Theorem}{Theorems}
\crefname{lem}{Lemma}{Lemmas}
\crefname{prop}{Proposition}{Propositions}
\crefname{eqnarray}{Problem}{Problems}

\begin{document}

\baselineskip 22pt \setcounter{page}{1}
\title{\bf A positive solution of the elliptic equation on a starshaped domain with boundary singularities
}
\author{{Zhi-Yun Tang, Xianhua Tang}\\
{\small \emph{School of Mathematics and Statistics, HNP-LAMA, Central South University,}}\\
{\small \emph{Changsha, Hunan 410083, People's Republic of China}}\\}
\date{}
\maketitle

{\bf Abstract}:
		We consider the elliptic equation with boundary singularities
		\begin{equation}
			\begin{cases}
				-\Delta u=-\lambda |x|^{-s_{1}}|u|^{p-2}u+|x|^{-s_{2}}|u|^{q-2}u &\text { in } \varOmega , \\
				u(x)=0 &\text { on } \partial \varOmega ,
			\end{cases}
		\end{equation}
		where $0\leq s_1 < s_2 < 2$, $2<p< 2^{*}(s_1)$, $q< 2^{*}(s_2)$. Which is the subcritical approximations of the Li-Lin's open problem proposed by Li and Lin (Arch Ration Mech Anal 203(3): 943-968, 2012). We find a positive solution which is a local minimum point of the energy functional on the Nehari manifold when $p>q>\frac{2-s_2}{2-s_1}p+\frac{2s_2-2s_1}{2-s_1}$. We also discuss the asymptotic behavior of the positive solution and find a new class of blow-up points by blowing up analysis. These blow-up points are on the boundary of the domain, which are not similar with the usual.



\section{Introduction}
Consider the elliptic equation with boundary singularities
		\begin{equation}
			\begin{cases}\label{eq1}
		-\Delta u=-\lambda |x|^{-s_{1}}|u|^{p-2}u+|x|^{-s_{2}}|u|^{q-2}u &\text { in } \varOmega , \\
								u(x)=0 &\text { on } \partial \varOmega ,
			\end{cases}
		\end{equation}
where $\varOmega  \subset \R $, $N\geq3$, $\lambda\in \RR$, $0\leq s_1 <s_2 < 2$, $2^{*}(s)=\frac{2(N-s)}{N-2}$ for $0 \leq s \leq 2, 2< p\leq 2^{*}(s_1), 2< q\leq 2^{*}(s_2).$

In \cite{LL2012}, Y. Li and C.-S. Lin asked an open question: Does the problem  have a positive solution when $\lambda>0$ and $p>q=2^{*}(s_2)$? The main difficult of this question is that one can not obtain the boundedness of the (PS) sequences for the energy functional. In 2015, G. Cerami, X. Zhong and W. Zou provided some partial positive answers in \cite{CZZ2015} by a perturbation approach and the monotonicity trick.

Recently in \cite{Tang2025}, we give a nonexistence result to the two-critical Li-Lin's open problem(the case that $\lambda>0$ and $2^{*}(s_1)=p>q=2^{*}(s_2)$) by the H\"{o}lder inequality, the Hardy inequality and Young inequality, and a positive solution for the Li-Lin's open problem(the case that $\lambda>0$ and $2^{*}(s_1)\geq p>q=2^{*}(s_2)$) by the method of sub-supersolutions. In \cite{Tang2024}, we obtain two positive solutions for small $\lambda$ by the least action principle and the mountain pass lemma, and no positive solutions for large $\lambda$ when  $2<q<\frac{2-s_2}{2-s_1}p+\frac{2s_2-2s_1}{2-s_1}$.

In this paper we consider the case that $q>\frac{2-s_2}{2-s_1}p+\frac{2s_2-2s_1}{2-s_1}$, obtain a positive solution by the Nehari manifold method. We also discuss the asymptotic behavior of the positive solution and find a new class of blow-up points by blowing up analysis. These blow-up points are on the boundary of the domain, which are not similar with the usual.  The main results are the following theorems.

\begin{thm}\label{thm1}
	Suppose that $\varOmega  \subset \R $ is a bounded starshaped domain about the origin. Assume that $\lambda>0$, $0\leq s_1 < s_2 < 2$, $2<p< 2^{*}(s_1)$, $q< 2^{*}(s_2)$ and $$p>q>\frac{2-s_2}{2-s_1}p+\frac{2s_2-2s_1}{2-s_1}.$$ Let
\begin{eqnarray}\label{01}
\frac{2^*(s_1)-p}{p(p-2)}>\frac{2^*(s_2)-q}{q(q-2)}.
\end{eqnarray}
Then \cref{eq1} has at least a positive solution $u_q \in C(\overline{\varOmega})$.
\end{thm}

\begin{rem} This question is the subcritical approximations of the Li-Lin's open problem.
\end{rem}

\begin{rem}
	For arbitrary $\lambda>0$, we obtain a positive solution. This is different from the case which has no positive solutions for large $\lambda$, where $2<q<\frac{2-s_2}{2-s_1}p+\frac{2s_2-2s_1}{2-s_1}$ as in \cite{Tang2024}.
\end{rem}

\begin{rem}
	We find a local minimum point of the Nehari manifold. In this case, there is no global minimum point.
\end{rem}

Not the same as subcritical situations, when $q$ is critical, we get a nonexistence result.
\begin{thm}\label{thm2}
	Suppose that $\varOmega  \subset \R $ is a bounded starshaped domain about the origin. Assume that $\lambda>0$, $0\leq s_1 < s_2 < 2$, $2^{*}(s_2) = q<p< 2^{*}(s_1)$.
	Then \cref{eq1} has no nonzero solution.
\end{thm}

Furthermore, we research the blowing up property for the positive solutions of \cref{thm1}.

\begin{thm}\label{thm3}
	Suppose that $\varOmega  \subset \R $ is a bounded starshaped domain about the origin. Assume that $0\leq s_1 < s_2 < 2$, $2^{*}(s_2) <p< 2^{*}(s_1)$ and $u_{q_j}$ is the positive solution of \cref{eq1} with $q=q_j$ by \cref{01}.
	Let $q_j \nearrow 2^*(s_2)$ and $\varepsilon_{j}=2^*(s_2)-q_j$. There exists a subsequence (still denoted by $\{q_j\}$) such that
	$$M_{q_j}=u_{q_{j}} ( x_{q_{j}} )=\operatorname* {m a x}\limits_{x \in \varOmega} u_{q_{j}}(x) \to\infty,\ x_{q_j} \to x_0\in \partial\Omega$$ as $j \to\infty$.
\end{thm}
\begin{rem}
	We can see that our blow up point $x_0\in \partial\Omega$, which is different from $x_0\in \Omega$ in \cite{H1991,R1989} due to the singularity of nonlinear term.
\end{rem}

\section{Preliminaries}
We  introduce the following work space
\begin{eqnarray*}
E=H^1_0{(\varOmega )}
\end{eqnarray*}
 with scalar product and norm given by
\begin{eqnarray*}
(u,v)=\int_{\varOmega }\nabla u\cdot\nabla vdx \quad\quad\quad \mbox{and}\quad\quad\quad {\|u\|}=(u,u)^{\frac{1}{2}}.
\end{eqnarray*}
It is well-known that the solutions of problem \eqref{eq1} are precisely the critical points of the energy  functional $I:H^1_0{(\varOmega )}\rightarrow \mathbb{R}$ defined by
\begin{eqnarray*}
I(u)= \frac{1}{2}{\|u\|}^{2}+\frac{\lambda}{p}\int_{\varOmega }|x|^{-s_{1}}|u|^p d x-\frac{1}{q}\int_{\varOmega }|x|^{-s_{2}}|u|^q d x.
\end{eqnarray*}
It is easy to see that $I$ is well-defined and $I\in C^2(H^1_0{(\varOmega )},\mathbb{R})$.  Then,   for any    $u,v,w\in H^1_0{(\varOmega )}$,
\begin{eqnarray*}
\langle I'(u),\ v \rangle =(u,v)+\lambda\int_{\varOmega }|x|^{-s_{1}}|u|^{p-2}uv dx- \int_{\varOmega }|x|^{-s_{2}}|u|^{q-2} u v d x,
\end{eqnarray*}
and
\begin{eqnarray*}
\langle I''(u)v,\ w\rangle =(v,w)+\lambda(p-1)\int_{\varOmega }|x|^{-s_{1}}|u|^{p-2}vw dx-(q-1) \int_{\varOmega }|x|^{-s_{2}}|u|^{q-2} v w d x.
\end{eqnarray*}
Let
\begin{eqnarray*}
\varphi(u)=\langle I'(u),\ u \rangle=\|u\|^2+\lambda\int_{\varOmega }|x|^{-s_{1}}|u|^pdx-\int_{\varOmega }|x|^{-s_{2}}|u|^qdx,
\end{eqnarray*}
and
\begin{eqnarray*}
\psi(u)=\langle I''(u)u,\ u\rangle=\|u\|^2+\lambda(p-1)\int_{\varOmega }|x|^{-s_{1}}|u|^pdx-(q-1)\int_{\varOmega }|x|^{-s_{2}}|u|^qdx.
\end{eqnarray*}
Define Nehari manifold
\begin{eqnarray*}
M=\left\{u\in H^1_0(\varOmega )\setminus \{0\}\left| \varphi(u)=0\right.\right\},
\end{eqnarray*}
and it's subsets
\begin{eqnarray*}
M^+=\left\{u\in M\left| \psi(u)<0\right.\right\},
\end{eqnarray*}
\begin{eqnarray*}
M^-=\left\{u\in M\left| \psi(u)>0\right.\right\},
\end{eqnarray*}
and
\begin{eqnarray*}
M^0=\left\{u\in M\left| \psi(u)=0\right.\right\}.
\end{eqnarray*}

\begin{prop}\label{prop1} Assume that $u\in H^1_0(\varOmega )\setminus \{0\}$. Then one has

(i) There exist $t^-_u>t^+_u>0$ such that $t^-_u u\in M^-$ and $t^+_u u\in M^+$ if
\begin{eqnarray*}
\|u\|^2-\frac{p-q}{p-2}\left(\frac{q-2}{p-2}\right)^\frac{q-2}{p-q}\left(\int_{\varOmega }|x|^{-s_{2}}|u|^qdx\right)^\frac{p-2}{p-q}\left(\lambda \int_{\varOmega }|x|^{-s_{1}}|u|^pdx\right)^{-\frac{q-2}{p-q}}<0.
\end{eqnarray*}

(ii) There exists $t^0_u>0$ such that $t^0_u u\in M^0$ if
\begin{eqnarray*}
\|u\|^2-\frac{p-q}{p-2}\left(\frac{q-2}{p-2}\right)^\frac{q-2}{p-q}\left(\int_{\varOmega }|x|^{-s_{2}}|u|^qdx\right)^\frac{p-2}{p-q}\left(\lambda \int_{\varOmega }|x|^{-s_{1}}|u|^pdx\right)^{-\frac{q-2}{p-q}}=0.
\end{eqnarray*}
\end{prop}

{\bf Proof\ \ \ \ } For $t>0$, define
\begin{eqnarray*}
h(t)=\frac{\varphi(tu)}{t^2}=\|u\|^2+\lambda t^{p-2} \int_{\varOmega }|x|^{-s_{1}}|u|^pdx-t^{q-2}\int_{\varOmega }|x|^{-s_{2}}|u|^qdx.
\end{eqnarray*}
We have
\begin{eqnarray*}
h'(t)=\lambda (p-2)t^{p-3}\int_{\varOmega }|x|^{-s_{1}}|u|^pdx-(q-2)t^{q-3}\int_{\varOmega }|x|^{-s_{2}}|u|^qdx,
\end{eqnarray*}
which implies that $h'(t_0)=0$ if and only if
\begin{eqnarray*}
t_0=\left(\frac{(q-2)\int_{\varOmega }|x|^{-s_{2}}|u|^qdx}{\lambda (p-2)\int_{\varOmega }|x|^{-s_{1}}|u|^pdx}\right)^\frac{1}{p-q}.
\end{eqnarray*}
Moreover, $h(t)$ is decreasing in $(0,t_0)$ and increasing in $(t_0,+\infty)$. Hence, $t_0$ is a minimum point of $h(t)$. Then we obtain
\begin{eqnarray*}
h(t_0)&=&\|u\|^2+\lambda \int_{\varOmega }|x|^{-s_{1}}|u|^pdx\left(\frac{(q-2)\int_{\varOmega }|x|^{-s_{2}}|u|^qdx}{\lambda (p-2)\int_{\varOmega }|x|^{-s_{1}}|u|^pdx}\right)^\frac{p-2}{p-q}\\&&-\int_{\varOmega }|x|^{-s_{2}}|u|^qdx\left(\frac{(q-2)\int_{\varOmega }|x|^{-s_{2}}|u|^qdx}{\lambda (p-2)\int_{\varOmega }|x|^{-s_{1}}|u|^pdx}\right)^\frac{q-2}{p-q}\\
&=&\|u\|^2-\frac{p-q}{p-2}\left(\frac{q-2}{p-2}\right)^\frac{q-2}{p-q}\left(\int_{\varOmega }|x|^{-s_{2}}|u|^qdx\right)^\frac{p-2}{p-q}\left(\lambda \int_{\varOmega }|x|^{-s_{1}}|u|^pdx\right)^{-\frac{q-2}{p-q}}.
\end{eqnarray*}
Let's prove the first one. In this case one has $h(t_0)<0$ and $h(0)=\|u\|^2>0$. Because $h(t)$ is decreasing in $(0,t_0)$, there exists a unique $t_u^+\in (0,t_0)$ such that $h(t_u^+)=0$, which implies that $\varphi(t_u^+ u)=0$, that is, $t_u^+ u \in M$. Furthermore we have
\begin{eqnarray*}
\psi((t_u^+)u)&=&(t_u^+)^2\|u\|^2+(p-1)(t_u^+)^p \lambda \int_{\varOmega }|x|^{-s_{1}}|u|^pdx-(q-1)(t_u^+)^q\int_{\varOmega }|x|^{-s_{2}}|u|^qdx\\
&=&(t_u^+)^2\left(h(t_u^+)+(p-2)(t_u^+)^{p-2}\lambda\int_{\varOmega }|x|^{-s_{1}}|u|^pdx-(q-2)(t_u^+)^{q-2}\int_{\varOmega }|x|^{-s_{2}}|u|^qdx\right)\\&=&(p-2)(t_u^+)^{p}\lambda\int_{\varOmega }|x|^{-s_{1}}|u|^pdx-(q-2)(t_u^+)^{q}\int_{\varOmega }|x|^{-s_{2}}|u|^qdx\\&=&\left((p-2)(t_u^+)^{p-q}\lambda\int_{\varOmega }|x|^{-s_{1}}|u|^pdx-(q-2)\int_{\varOmega }|x|^{-s_{2}}|u|^qdx\right) (t_u^+)^{q}\\
&<&\left((p-2)t_0^{p-q}\lambda\int_{\varOmega }|x|^{-s_{1}}|u|^pdx-(q-2)\int_{\varOmega }|x|^{-s_{2}}|u|^qdx\right) (t_u^+)^{q}\\
&=&0,
\end{eqnarray*}
which implies that $t^+_u u\in M^+$. Similarly, one can prove there exists a unique $t_u^-\in (t_0,+\infty)$ such that $t^-_u u\in M^-$.

Next we consider the second case. In this case one has $h(t_0)=0$. Let $t_u^0=t_0$, then $\varphi(t_u^0 u)=0$, that is, $t_u^0 u \in M$. Moreover we have
\begin{eqnarray*}
	\psi(t_u^0u)&=&(t_u^0)^2\|u\|^2+(p-1)(t_u^0)^p \lambda \int_{\varOmega }|x|^{-s_{1}}|u|^pdx-(q-1)(t_u^0)^q\int_{\varOmega }|x|^{-s_{2}}|u|^qdx\\
	&=&(t_u^0)^2\left(h(t_u^0)+(p-2)(t_u^0)^{p-2}\lambda\int_{\varOmega }|x|^{-s_{1}}|u|^pdx-(q-2)(t_u^0)^{q-2}\int_{\varOmega }|x|^{-s_{2}}|u|^qdx\right)\\
	&=&(p-2)(t_u^0)^{p}\lambda\int_{\varOmega }|x|^{-s_{1}}|u|^pdx-(q-2)(t_u^0)^{q}\int_{\varOmega }|x|^{-s_{2}}|u|^qdx\\
	&=&\left((p-2)(t_0)^{p-q}\lambda\int_{\varOmega }|x|^{-s_{1}}|u|^pdx-(q-2)\int_{\varOmega }|x|^{-s_{2}}|u|^qdx\right) (t_u^0)^{q}\\
	&=&0,
\end{eqnarray*}
which implies that $t^0_u u\in M^0$. We complete the proof.
$\hfill\Box$

\begin{prop}\label{p2}
	Define
	$$
	m^+=\inf\{I(u)\mid u\in M^+\}, \ \ m^0=\inf\{I(u)\mid u\in M^0\}.
	$$
	Then we have \label{m^+} $m^+>0$ if $M^+\not=\phi$ and $m^0>0$ if $M^0\not=\phi$.
\end{prop}
{\bf Proof\ \ \ \ } For $u\in M^+$, we have
\begin{align}
&\|u\|^2+\lambda\int_{\varOmega }|x|^{-s_{1}}|u|^pdx-\int_{\varOmega }|x|^{-s_{2}}|u|^qdx=\varphi(u)=0,\label{2.0}\\
&\|u\|^2+\lambda(p-1)\int_{\varOmega }|x|^{-s_{1}}|u|^pdx-(q-1)\int_{\varOmega }|x|^{-s_{2}}|u|^qdx=\psi(u)<0\nonumber,
\end{align}
which implies that
\begin{eqnarray*}
\lambda\int_{\varOmega }|x|^{-s_{1}}|u|^pdx&=&\frac{q-2}{p-q}\|u\|^2+\frac{\psi(u)}{p-q},\\
\int_{\varOmega }|x|^{-s_{2}}|u|^qdx&=&\frac{p-2}{p-q}\|u\|^2+\frac{\psi(u)}{p-q}.
\end{eqnarray*}
The one obtains
\begin{eqnarray*}
I(u)&=&\frac{1}{2}\|u\|^2+\frac{1}{p}\lambda\int_{\varOmega }|x|^{-s_{1}}|u|^pdx-\frac{1}{q}\int_{\varOmega }|x|^{-s_{2}}|u|^qdx\\
&=&\frac{(p-2)(q-2)}{2pq}\|u\|^2-\frac{\psi(u)}{pq}\\
&>&\frac{(p-2)(q-2)}{2pq}\|u\|^2.
\end{eqnarray*}
It follows from (\ref{2.0}) and Hardy-Sobolev inequality that
\begin{eqnarray*}
\|u\|^2&\leq&\|u\|^2+\lambda\int_{\varOmega }|x|^{-s_{1}}|u|^pdx\\
&=&\int_{\varOmega }|x|^{-s_{2}}|u|^qdx\\
&\leq&C\|u\|^q
\end{eqnarray*}
for all $u\in M$ and some constant $C>0$, which implies that $\|u\|\geq C^{-\frac{1}{q-2}}$. Hence,
\begin{equation}\label{eq2.2}
	I(u)>\frac{(p-2)(q-2)}{2pq}\|u\|^2\geq \frac{(p-2)(q-2)}{2pq}C^{-\frac{2}{q-2}}>0
\end{equation}
for all $u \in M^+$. Then $m^+>0.$ Similarly, we have $m^0>0$.
$\hfill\Box$

\begin{defn}
We say that $\varOmega  \subset \R$ is a starshaped domain about the origin, if $sx \in \varOmega $ for all $s \in (0,1) $ and $x \in \varOmega $.

	For $0<r\leq 1$ and $u\in H^1_0(\varOmega )$, define
	$$
	u_{r}(x)=u(r^{-1}x)$$
	for all $x \in \varOmega $. Then $u_r\in H^1_0(\varOmega )$ for $0<r\leq 1$. Moreover, we have
	\begin{eqnarray*}
		\int_{\varOmega }|\nabla u_r(x)|^2dx&=&r^{N-2}\int_{\varOmega }|\nabla u(x)|^2dx,
	\end{eqnarray*}
	\begin{eqnarray*}
		\int_{\varOmega }|x|^{-s_1}|u_{r}(x)|^{p}dx&=& r^{N-s_1}\int_{\varOmega }|x|^{-s_1}|u(x)|^{p}dx,
	\end{eqnarray*}
	and
	\begin{eqnarray*}
		\int_{\varOmega }|x|^{-s_2}|u_{r}(x)|^{q}dx
		&=&r^{N-s_2}\int_{\varOmega }|x|^{-s_2}|u(x)|^qdx.
	\end{eqnarray*}
\end{defn}

\begin{prop}\label{M+}
Suppose that $\varOmega  \subset \R$ is a bounded starshaped domain about the origin. Assume that $0\leq s_1 < s_2 < 2$, $2<p\leq 2^{*}(s_1)$, $q\leq 2^{*}(s_2)$ and $p>q>\frac{2-s_2}{2-s_1}p+\frac{2s_2-2s_1}{2-s_1}$. Then we have $M^+\not=\emptyset$.
\end{prop}

{\bf Proof}\ \ \ \
Then for every $u\in H^1_0(\varOmega )\setminus \{0\}$ and $0<r\leq 1$, we have
\begin{eqnarray*}
&&\|u_r\|^2-\frac{p-q}{p-2}\left(\frac{q-2}{p-2}\right)^\frac{q-2}{p-q}\left(\int_{\varOmega }|x|^{-s_{2}}|u_r|^qdx\right)^\frac{p-2}{p-q}\left(\lambda \int_{\varOmega }|x|^{-s_{1}}|u_r|^pdx\right)^{-\frac{q-2}{p-q}}\\
&=&r^{N-2}\|u\|^2-\frac{p-q}{p-2}\left(\frac{q-2}{p-2}\right)^\frac{q-2}{p-q}\left(r^{N-s_2}\int_{\varOmega }|x|^{-s_{2}}|u|^qdx\right)^\frac{p-2}{p-q}\left(r^{N-s_1}\lambda \int_{\varOmega }|x|^{-s_{1}}|u|^pdx\right)^{-\frac{q-2}{p-q}}\\
&=&r^{N-2}\left\{\|u\|^2-r^{\alpha}\frac{p-q}{p-2}\left(\frac{q-2}{p-2}\right)^\frac{q-2}{p-q}\left(\int_{\varOmega }|x|^{-s_{2}}|u|^qdx\right)^\frac{p-2}{p-q}\left(\lambda \int_{\varOmega }|x|^{-s_{1}}|u|^pdx\right)^{-\frac{q-2}{p-q}}\right\}\\
&\leq&0
\end{eqnarray*}
for $\alpha\stackrel{\bigtriangleup}{=}\frac{(p-2)(2-s_2)-(q-2)(2-s_1)}{p-q}<0$ and $0<r\leq 1$ small enough, which implies that exists unique $t^+_{u_r}>0$ such that $t^+_{u_r} u_r\in M^+$ by Proposition \ref{prop1}.
Then we have $M^+\not=\emptyset$. $\hfill\Box$

\begin{prop}\label{p5}
Suppose that $\varOmega  \subset \R$ is a bounded star-shaped domain with star-concenter 0. Assume that $0\leq s_1 < s_2 < 2$, $2<p< 2^{*}(s_1)$, $q< 2^{*}(s_2)$ and $p>q>\frac{2-s_2}{2-s_1}p+\frac{2s_2-2s_1}{2-s_1}$. Let
\begin{eqnarray*}
\frac{2^*(s_1)-p}{p(p-2)}>\frac{2^*(s_2)-q}{q(q-2)}.
\end{eqnarray*}
Then we have $m^+<m^0$ if $M^0\not=\phi$.
\end{prop}

{\bf Proof}\ \ \ \
We choose $u_n\in M^0$ such that $I(u_n)\to m^0$ as $n\to\infty$. Then one has
\begin{eqnarray}
\|u_n\|^2+\lambda\int_{\varOmega }|x|^{-s_{1}}|u_n|^pdx-\int_{\varOmega }|x|^{-s_{2}}|u_n|^qdx&=&0,\label{2.1}\\
\|u_n\|^2+\lambda(p-1)\int_{\varOmega }|x|^{-s_{1}}|u_n|^pdx-(q-1)\int_{\varOmega }|x|^{-s_{2}}|u_n|^qdx&=&0,\nonumber
\end{eqnarray}
which implies that
\begin{eqnarray}
\lambda\int_{\varOmega }|x|^{-s_{1}}|u_n|^pdx&=&\frac{q-2}{p-q}\|u_n\|^2,\label{2.2}\\
\int_{\varOmega }|x|^{-s_{2}}|u_n|^qdx&=&\frac{p-2}{p-q}\|u_n\|^2\label{2.3}
\end{eqnarray}
and \begin{eqnarray}\label{pq}
\lambda\int_{\varOmega }|x|^{-s_{1}}|u_n|^pdx&=&
\frac{q-2}{p-2}\int_{\varOmega }|x|^{-s_{2}}|u_n|^qdx.
\end{eqnarray}
Hence, we obtain
\begin{eqnarray}\label{m^0}
\frac{(p-q)(q-2)}{2pq}\int_{\varOmega }|x|^{-s_{2}}|u_n|^qdx=\frac{(p-2)(q-2)}{2pq}\|u_n\|^2
=I(u_n)\to m^0.
\end{eqnarray}
It follows that $\{u_n\}$ is bounded. Going if necessary to a subsequence, we can assume that
\begin{eqnarray*}
&&u_n\rightharpoonup u \ \ \ \ \  \ \ \ \ \ \  \ \ \  \ \text{in}\ H^1_0{(\varOmega)},\\
&&u_n\rightarrow u \ \ \ \ \ \ \ \ \  \ \ \ \ \ \ \text{in}\ L^{r}( \varOmega)\ (r\in [2,2^*)),\\
&&u_n(x)\rightarrow u(x) \ \ \ \ \  \ \ a.e. \ \text{in}\  \varOmega,\\
&&\int_{\varOmega }|x|^{-s_{1}}|u_n|^pdx\rightarrow\int_{\varOmega }|x|^{-s_{1}}|u|^pdx,\\
&&\int_{\varOmega }|x|^{-s_{2}}|u_n|^qdx\rightarrow\int_{\varOmega }|x|^{-s_{2}}|u|^qdx
\end{eqnarray*}
as $n\rightarrow\infty$. Moreover, $\lim\limits_{n\to\infty}\|u_n\|\geq\|u\|$. We consider two cases that $\lim\limits_{n\to\infty}\|u_n\|>\|u\|$ and $\lim\limits_{n\to\infty}\|u_n\|=\|u\|$ respectively.

In the case that $\lim\limits_{n\to\infty}\|u_n\|>\|u\|$. By (\ref{pq}) one has
\begin{eqnarray}\label{eq3}
\lambda\int_{\varOmega }|x|^{-s_{1}}|u|^pdx=
\frac{q-2}{p-2}\int_{\varOmega }|x|^{-s_{2}}|u|^qdx.
\end{eqnarray}
It follows from (\ref{m^0}) that
\begin{eqnarray}\label{eq2}
	m^0=\frac{(p-q)(q-2)}{2pq}\int_{\varOmega }|x|^{-s_{2}}|u|^qdx.
\end{eqnarray}
Hence, $u\not=0$ because that $m^0>0$ by Proposition \ref{p2}.
From (\ref{2.2}) and (\ref{2.3}) we obtain
\begin{eqnarray*}
\|u_n\|^2-\frac{p-q}{p-2}\left(\frac{q-2}{p-2}\right)^\frac{q-2}{p-q}\left(\int_{\varOmega }|x|^{-s_{2}}|u_n|^qdx\right)^\frac{p-2}{p-q}\left(\lambda \int_{\varOmega }|x|^{-s_{1}}|u_n|^pdx\right)^{-\frac{q-2}{p-q}}=0,
\end{eqnarray*}
which implies that
\begin{eqnarray*}
\|u\|^2-\frac{p-q}{p-2}\left(\frac{q-2}{p-2}\right)^\frac{q-2}{p-q}\left(\int_{\varOmega }|x|^{-s_{2}}|u|^qdx\right)^\frac{p-2}{p-q}\left(\lambda \int_{\varOmega }|x|^{-s_{1}}|u|^pdx\right)^{-\frac{q-2}{p-q}}<0.
\end{eqnarray*}
There exists unique $t_u^+>0$ such that $t_u^+ u\in M^+$ by  Proposition \ref{prop1}. Due to (\ref{2.1}), we have
$$\|u\|^2+\lambda \int_{\varOmega }|x|^{-s_{1}}|u|^pdx-\int_{\varOmega }|x|^{-s_{2}}|u|^qdx<0.$$
By $t_u^+u \in M^+ \subset M$, one has
$$\varphi(t_u^+u)=(t_u^+)^2\|u\|^2+(t_u^+)^p\lambda \int_{\varOmega }|x|^{-s_{1}}|u|^pdx-(t_u^+)^q\int_{\varOmega }|x|^{-s_{2}}|u|^qdx=0,$$
which implies that $t_u^+\not=1$. Furthermore, we obtain
\begin{eqnarray*}
I(t_u^+u)
&=&\frac{1}{2}(t_u^+)^2\|u\|^2+\frac{1}{p}\lambda (t_u^+)^p\int_{\varOmega }|x|^{-s_{1}}|u|^pdx-\frac{1}{q}(t_u^+)^q\int_{\varOmega }|x|^{-s_{2}}|u|^qdx-\dfrac{1}{2}\varphi(t_u^+u)\\
&=&-\frac{p-2}{2p}\lambda (t_u^+)^p\int_{\varOmega }|x|^{-s_{1}}|u|^pdx+\frac{q-2}{2q}(t_u^+)^q\int_{\varOmega }|x|^{-s_{2}}|u|^qdx\\
&=&-\frac{q-2}{2p} (t_u^+)^p\int_{\varOmega }|x|^{-s_{2}}|u|^qdx+\frac{q-2}{2q}(t_u^+)^q\int_{\varOmega }|x|^{-s_{2}}|u|^qdx
\end{eqnarray*}
by (\ref{eq3}). Consider the function
\begin{eqnarray*}h(t)&\stackrel{\triangle}{=}&-\frac{q-2}{2p} t^p\int_{\varOmega }|x|^{-s_{2}}|u|^qdx+\frac{q-2}{2q}t^q\int_{\varOmega }|x|^{-s_{2}}|u|^qdx,
\end{eqnarray*}
we have
\begin{eqnarray*}h'(t)&=&-\frac{q-2}{2} t^{q-1}(t^{p-q}-1)\int_{\varOmega }|x|^{-s_{2}}|u|^qdx.
\end{eqnarray*}
Hence, $h$ is strictly increasing on $[0,1]$ and strictly decreasing on $[1,+\infty)$, which implies that
$$
m^+\leq I(t_u^+u)=h(t_u^+)<h(1)=m^0
$$
by (\ref{eq2}).

In the case that $\lim\limits_{n\to\infty}\|u_n\|=\|u\|$, we have $u_n \to u$, $u\in M^0$ and $I(u)=m^0$.
Assume that $0<r\leq 1$ and $t>0$. Then $tu_r\in M$ if and only if
$$t^2\|u_r\|^2+\lambda t^p\int_{\varOmega }|x|^{-s_{1}}|u_r|^pdx-t^q\int_{\varOmega }|x|^{-s_{2}}|u_r|^qdx=0. $$
Note that
\begin{eqnarray*}
&&t^2\|u_r\|^2+\lambda t^p\int_{\varOmega }|x|^{-s_{1}}|u_r|^pdx-t^q\int_{\varOmega }|x|^{-s_{2}}|u_r|^qdx\\
&=&t^2r^{N-2}\|u\|^2+\lambda t^pr^{N-s_1}\int_{\varOmega }|x|^{-s_{1}}|u|^pdx-t^qr^{N-s_2}\int_{\varOmega }|x|^{-s_{2}}|u|^qdx\\
&=&\left(t^2r^{N-2}+\frac{q-2}{p-q}t^{p}r^{N-s_1}-\frac{p-2}{p-q}t^{q}r^{N-s_2}\right)\|u\|^2\\
&=&\left(1+\frac{q-2}{p-q}t^{p-2}r^{2-s_1}-\frac{p-2}{p-q}t^{q-2}r^{2-s_2}\right)t^2r^{N-2}\|u\|^2.
\end{eqnarray*}
Let
\begin{eqnarray*}
h(t,r)&=&1+\frac{q-2}{p-q}t^{p-2}r^{2-s_1}-\frac{p-2}{p-q}t^{q-2}r^{2-s_2}.
\end{eqnarray*}
Due to $h(1,1)=0, \frac{\partial h}{\partial r}(1,1)=\frac{(q-2)(2-s_1)-(p-2)(2-s_2)}{p-q}>0, h\in C^2(\RR^+\times \RR^+,\RR)$, it follows from the implicit function theorem that there exists a constant $\delta_0>0$ and a $C^2$ function $r=r(t)$ from $(1-\delta_0, 1+\delta_0)$ to $(1-\delta_0, 1+\delta_0)$ such that $r(1)=1$ and $h(t,r(t))=0$, which implies that $tu_r\in M$. Furthermore, one has
\begin{eqnarray*}
r'(1)&=&\left.-\frac{h_t}{h_r} \right|_{t=1}=0,\\
r''(1)&=&\left. \frac{2h_th_rh_{tr}-h_r^2h_{tt}-h_t^2h_{rr}}{h_r^3} \right|_{t=1} \\
&=&\left. \frac{-h_{tt}}{h_r}\right|_{t=1}\\
&=&\frac{-(p-2)(q-2)(p-q)}{(q-2)(2-s_1)-(p-2)(2-s_2)}<0.
\end{eqnarray*}
\begin{eqnarray*}
	\psi(tu_r)&=&t^2\|u_r\|^2+(p-1)t^p \lambda \int_{\varOmega }|x|^{-s_{1}}|u_r|^pdx-(q-1)t^q\int_{\varOmega }|x|^{-s_{2}}|u_r|^qdx\\
	&=&t^2r^{N-2}\|u\|^2+(p-1)t^pr^{N-s_1} \lambda \int_{\varOmega }|x|^{-s_{1}}|u|^pdx-(q-1)t^qr^{N-s_2}\int_{\varOmega }|x|^{-s_{2}}|u|^qdx\\
	&=&\left(t^2r^{N-2}+\frac{(p-1)(q-2)}{p-q}t^{p}r^{N-s_1}-\frac{(q-1)(p-2)}{p-q}t^{q}r^{N-s_2}\right)\|u\|^2\\
	&=&\left(1+\frac{(p-1)(q-2)}{p-q}t^{p-2}r^{2-s_1}-\frac{(q-1)(p-2)}{p-q}t^{q-2}r^{2-s_2}\right)t^2r^{N-2}\|u\|^2\\		&=&\left(1+\frac{(p-1)(q-2)}{p-q}t^{p-2}r^{2-s_1}-\frac{(q-1)(p-2)}{p-q}t^{q-2}r^{2-s_2}-h(t,r)\right)t^2r^{N-2}\|u\|^2\\
	&=&\left(\frac{(p-2)(q-2)}{p-q}t^{p-2}r^{2-s_1}-\frac{(p-2)(q-2)}{p-q}t^{q-2}r^{2-s_2}\right)t^2r^{N-2}\|u\|^2\\
	&=&\left(t^{p}r^{N-s_1}-t^{q}r^{N-s_2}\right)\frac{(p-2)(q-2)}{p-q}\|u\|^2
\end{eqnarray*}
by (\ref{2.2}), (\ref{2.3}) and $h(t,r)=0$ .

Define
\begin{eqnarray*}
	f(t)&=&t^{p}r^{N-s_1}-t^{q}r^{N-s_2}.
\end{eqnarray*}
Then one has
\begin{eqnarray*}
	f'(t)&=&pt^{p-1}r^{N-s_1}-qt^{q-1}r^{N-s_2}+(N-s_1)t^{p}r^{N-1-s_1}r'-(N-s_2)t^{q}r^{N-1-s_2}r'.
\end{eqnarray*}
Hence, we have
\begin{eqnarray*}
	f'(1)&=&p-q>0.
\end{eqnarray*}
Then there exists a constant $\delta_1 \in (0,\delta_0)$ such that $f'(t)\geq \frac{1}{2}(p-q)>0$ for $t \in (1-\delta_1,1]$, which implies that $f(t)<f(1)=0$ for $t \in (1-\delta_1,1)$. Thus, $tu_r \in M^+$ for $t \in (1-\delta_1,1)$. Moreover, we obtain
\begin{eqnarray*}
I(tu_r)
&=&\frac{1}{2}t^2\|u_r\|^2+\frac{1}{p}\lambda t^p\int_{\varOmega }|x|^{-s_{1}}|u_r|^pdx-\frac{1}{q}t^q\int_{\varOmega }|x|^{-s_{2}}|u_r|^qdx\\
&=&\frac{1}{2}t^2r^{N-2}\|u\|^2+\frac{1}{p}\lambda t^pr^{N-s_1}\int_{\varOmega }|x|^{-s_{1}}|u|^pdx-\frac{1}{q}t^qr^{N-s_2}\int_{\varOmega }|x|^{-s_{2}}|u|^qdx\\
&=&\left(\frac{1}{2}t^2r^{N-2}+\frac{q-2}{p(p-q)}t^{p}r^{N-s_1}-\frac{p-2}{q(p-q)}t^{q}r^{N-s_2}\right)\|u\|^2\\
&=&\left(\frac{1}{2}+\frac{q-2}{p(p-q)}t^{p-2}r^{2-s_1}-\frac{p-2}{q(p-q)}t^{q-2}r^{2-s_2}-\dfrac{1}{2}h(t,r)\right)t^2r^{N-2}\|u\|^2\\
&=&\left(-\frac{(p-2)(q-2)}{2p(p-q)}t^{p-2}r^{2-s_1}+\frac{(p-2)(q-2)}{2q(p-q)}t^{q-2}r^{2-s_2}\right)t^2r^{N-2}\|u\|^2\\
&=&\left(-qt^{p}r^{N-s_1}+pt^{q}r^{N-s_2}\right)\frac{(p-2)(q-2)}{2pq(p-q)}\|u\|^2
\end{eqnarray*}
by (\ref{2.2}), (\ref{2.3}) and $h(t,r)=0$. Define
\begin{eqnarray*}
g(t)&=&-qt^{p}r^{N-s_1}+pt^{q}r^{N-s_2}.
\end{eqnarray*}
Then one has
\begin{eqnarray*}
g'(t)&=&-pqt^{p-1}r^{N-s_1}+pqt^{q-1}r^{N-s_2}-q(N-s_1)t^{p}r^{N-1-s_1}r'+p(N-s_2)t^{q}r^{N-1-s_2}r',\\
g''(t)&=&-q[p(p-1)t^{p-2}r^{N-s_1}+2p(N-s_1)t^{p-1}r^{N-1-s_1}r'\\
&&+(N-s_1)(N-1-s_1)t^{p}r^{N-2-s_1}r'^2+(N-s_1)t^{p}r^{N-1-s_1}r'']\\
&&+p[q(q-1)t^{q-2}r^{N-s_2}+2q(N-s_2)t^{q-1}r^{N-1-s_2}r'\\
&&+(N-s_2)(N-1-s_2)t^{q}r^{N-2-s_2}r'^2+(N-s_2)t^{q}r^{N-1-s_2}r''].
\end{eqnarray*}
Hence, we have
\begin{eqnarray*}
g'(1)&=&0,\\
g''(1)&=&-q[p(p-1)+(N-s_1)r''(1)]+p[q(q-1)+(N-s_2)r''(1)]\\
&=&-qp(p-q)+[p(N-s_2)-q(N-s_1)]r''(1)\\
&=&-qp(p-q)-[p(N-s_2)-q(N-s_1)]\frac{(p-2)(q-2)(p-q)}{(q-2)(2-s_1)-(p-2)(2-s_2)}\\
&=&\frac{p-q}{(q-2)(2-s_1)-(p-2)(2-s_2)}\\&&\left\lbrace  -qp\left[ (q-2)(2-s_1)-(p-2)(2-s_2)\right] -[p(N-s_2)-q(N-s_1)](p-2)(q-2) \right\rbrace\\
&=&\frac{pq(p-q)(q-2)(p-2)}{(q-2)(2-s_1)-(p-2)(2-s_2)}\left\lbrace  -\frac{2-s_1}{p-2}+\frac{2-s_2}{q-2} -\frac{N-s_2}{q}+\frac{N-s_1}{p} \right\rbrace\\
&=&\frac{pq(p-q)(q-2)(p-2)}{(q-2)(2-s_1)-(p-2)(2-s_2)}\\&&\left\lbrace  \frac{(N-s_1)(p-2)-p(2-s_1)}{p(p-2)}-\frac{(N-s_2)(q-2)-q(2-s_2)}{q(q-2)} \right\rbrace\\
&=&\frac{pq(p-q)(q-2)(p-2)}{(q-2)(2-s_1)-(p-2)(2-s_2)}\left\lbrace  \frac{p(N-2)-2(N-s_1)}{p(p-2)}-\frac{q(N-2)-2(N-s_2)}{q(q-2)} \right\rbrace\\
&=&-\frac{pq(p-q)(q-2)(p-2)(N-2)}{(q-2)(2-s_1)-(p-2)(2-s_2)}\left\lbrace  \frac{2^{*}(s_1)-p}{p(p-2)}-\frac{2^{*}(s_2)-q}{q(q-2)} \right\rbrace\\
&<&0
\end{eqnarray*}
from (\ref{01}). Then there exists a constant $\delta \in (0,\delta_1)$ such that $g''(t)\leq \frac{1}{2}g''(1)<0$ for $t \in (1-\delta,1]$, which implies that $g(t)<g(1)=p-q$ for $t \in (1-\delta,1)$. Thus, we have
\begin{eqnarray*}
I(tu_r)&=&\left(-qt^{p}r^{N-s_1}+pt^{q}r^{N-s_2}\right)\frac{(p-2)(q-2)}{2pq(p-q)}\|u\|^2\\
&=&g(t)\frac{(p-2)(q-2)}{2pq(p-q)}\|u\|^2\\
&<&g(1)\frac{(p-2)(q-2)}{2pq(p-q)}\|u\|^2\\
&=&\frac{(p-2)(q-2)}{2pq}\|u\|^2\\
&=&I(u)
\end{eqnarray*}
for $t \in (1-\delta,1)$ by (\ref{m^0}). Hence, we obtain
\begin{eqnarray*}
	m^+\leq I(tu_r)<I(u)=m_0.
\end{eqnarray*}
Then $m^+<m^0$.$\hfill\Box$
\begin{prop}\label{p6}
	\ \ \ \ Let $m^+_q=m^+$, $M^+_q=M^+$, $I_q(u)=I(u)$ and $\psi_q(u)=\psi(u)$ with $2<q\leq 2^*(s_2)$. Because that $u_{q_j}$ is the positive solution of \cref{eq1} by \cref{01}, we have
	$$I_{q_j}'(u_{q_j})=0,~I_{q_j}(u_{q_j})=m_{q_j}^+,~\psi_{q_j}(u_{q_j})\leq 0.$$
	Then $\limsup\limits_{j \to \infty} m_{q_j}^+\leq m_{2^*(s_2)}^+$.
\end{prop}

{\bf Proof}\ \ \ \ From the Proposition \ref{p2} and Proposition \ref{M+}, we get $M_{2^*(s_2)}\neq \emptyset$  and $0<m_{2^*(s_2)}^+<\infty$. By the definition of $m_{2^*(s_2)}^+$, for every $\varepsilon>0$, there exists $u \in M_{2^*(s_2)}^+$ such that
$$I_{2^*(s_2)}(u) < m_{2^*(s_2)}^+ +\varepsilon.$$
Define $h(t, q)=\varphi_q(tu)$. Then $h(1,2^*(s_2))=0$, $$h'_t (1,2^*(s_2))=\langle \varphi'_{2^*(s_2)}(u), u \rangle =\psi_{2^*(s_2)}(u)+\varphi_{2^*(s_2)}(u)=\psi_{2^*(s_2)}(u)<0.$$
Due to $h\in C^2(\RR^+\times \RR^+,\RR)$, it follows from the implicit function theorem that there exists a constant $\delta_0>0$ and a $C^2$ function $t=t(q)$ from $(2^*(s_2)-\delta_0, 2^*(s_2)+\delta_0)$ to $(1-\delta_0, 1+\delta_0)$ such that $t(2^*(s_2))=1$ and $h(t(q),q)=0$, which implies that $t(q)u\in M_q$.

Define $f(q)=\psi(t(q)u)$, then $f(2^*(s_2))=\psi(u)<0$ due to $u \in M_{2^*(s_2)}^+$. It is easy to know $f(q)$ is continuous, there exists $\delta_1 \in (0,\delta_0)$ such that $f(q)<0$ for $q \in (2^*(s_2)-\delta_1, 2^*(s_2))$, that is, $t(q)u \in M^+_q$. Hence, we have
$$m_q^+ \leq I_q(t(q)u) \leq I_{2^*(s_2)}(t({2^*(s_2)})u)+\varepsilon=I_{2^*(s_2)}(u)+\varepsilon \leq m_{2^*(s_2)}^+ + 2\varepsilon, $$
which implies that
$$\limsup\limits_{q \to {2^*(s_2)}} m_q^+ \leq m_{2^*(s_2)}^+ + 2\varepsilon. $$
By the arbitrariness of $\varepsilon$, one has $\limsup\limits_{j \to \infty} m_{q_j}^+\leq m_{2^*(s_2)}^+$.$\hfill\Box$

{\bf Proof of Theorem \ref{thm1}}\ \ \ \ Choose $u_n\in M^+$ such that $I(u_n)\to m^+$ as $n\to\infty$. Without loss of generality, we assume $u_n\geq 0$ because $|u|\in M^+$ when $u \in M^+ $ and $I(|u|)=I(u)$ when $u \in \HO $. Then one has
\begin{eqnarray}
\varphi(u_n)=\|u_n\|^2+\lambda\int_{\varOmega }|x|^{-s_{1}}|u_n|^pdx-\int_{\varOmega }|x|^{-s_{2}}|u_n|^qdx&=&0,\\
\psi(u_n)=\|u_n\|^2+(p-1)\lambda\int_{\varOmega }|x|^{-s_{1}}|u_n|^pdx-(q-1)\int_{\varOmega }|x|^{-s_{2}}|u_n|^qdx&<&0\nonumber,
\end{eqnarray}
which implies that
\begin{eqnarray*}
\lambda\int_{\varOmega }|x|^{-s_{1}}|u_n|^pdx&=&\frac{q-2}{p-q}\|u_n\|^2+\frac{\psi(u_n)}{p-q},\\
\int_{\varOmega }|x|^{-s_{2}}|u_n|^qdx&=&\frac{p-2}{p-q}\|u_n\|^2+\frac{\psi(u_n)}{p-q}.
\end{eqnarray*}
The one obtains
\begin{eqnarray*}
I(u_n)&=&\frac{1}{2}\|u_n\|^2+\frac{1}{p}\lambda\int_{\varOmega }|x|^{-s_{1}}|u_n|^pdx-\frac{1}{q}\int_{\varOmega }|x|^{-s_{2}}|u_n|^qdx\\
&=&\frac{(p-2)(q-2)}{2pq}\|u_n\|^2-\frac{\psi(u_n)}{pq}\\
&\geq&\frac{(p-2)(q-2)}{2pq}\|u_n\|^2.
\end{eqnarray*}
It follows that $\{u_n\}$ is bounded. Going if necessary to a subsequence, we can assume that
\begin{eqnarray*}
&&u_n\rightharpoonup u \ \ \ \ \  \ \ \ \ \ \  \ \ \  \ \text{in}\ H^1_0{(\varOmega)},\\
&&u_n\rightarrow u \ \ \ \ \ \ \ \ \  \ \ \ \ \ \ \text{in}\ L^{q}( \varOmega)\ (q\in [2,2^*)),\\
&&u_n(x)\rightarrow u(x) \ \ \ \ \  \ \ a.e. \ \text{in}\  \varOmega,\\
&&\int_{\varOmega }|x|^{-s_{1}}|u_n|^pdx\rightarrow\int_{\varOmega }|x|^{-s_{1}}|u|^pdx,\\
&&\int_{\varOmega }|x|^{-s_{2}}|u_n|^qdx\rightarrow\int_{\varOmega }|x|^{-s_{2}}|u|^qdx
\end{eqnarray*}
as $n\rightarrow\infty$. Moreover, $u\geq 0$ and  $\lim\limits_{n\to\infty}\|u_n\|\geq\|u\|$.

Assume that $\lim\limits_{n\to\infty}\|u_n\|>\|u\|$. Let $\psi(u_n)\to -m$. Then $m\geq 0$. It follows that
\begin{eqnarray*}
-m&=&\lim_{n\to\infty}(\psi(u_n)-\varphi(u_n))=(p-2)\lambda \int_{\varOmega }|x|^{-s_{1}}|u|^pdx-(q-2)\int_{\varOmega }|x|^{-s_{2}}|u|^qdx,\\
m^+&=&\lim_{n\to\infty}\left( I(u_n)-\dfrac{\varphi(u_n)}{2}\right) =-\frac{p-2}{2p}\lambda \int_{\varOmega }|x|^{-s_{1}}|u|^pdx+\frac{q-2}{2q}\int_{\varOmega }|x|^{-s_{2}}|u|^qdx,
\end{eqnarray*}
which implies that
\begin{eqnarray}
\frac{(p-2)(p-q)}{2pq}\lambda \int_{\varOmega }|x|^{-s_{1}}|u|^pdx&=&m^+-\frac{m}{2q},\label{10}\\
\frac{(q-2)(p-q)}{2pq} \int_{\varOmega }|x|^{-s_{2}}|u|^qdx&=&m^+-\frac{m}{2p}\label{11}.
\end{eqnarray}
Let
\begin{eqnarray*}
g(t)
&=&\|u\|^2+\lambda t^{p-2}\int_{\varOmega }|x|^{-s_{1}}|u|^pdx-t^{q-2}\int_{\varOmega }|x|^{-s_{2}}|u|^qdx.
\end{eqnarray*}
Due to (\ref{2.1}), we have
$$\|u\|^2+\lambda \int_{\varOmega }|x|^{-s_{1}}|u|^pdx-\int_{\varOmega }|x|^{-s_{2}}|u|^qdx<0,$$ which implies that $g(1)<0$ and $u\neq0$. Note that $g(0)>0$. Hence, there exists $t_u^+\in (0, 1)$ such that $g(t_u^+)=0$, that is,
$$\varphi(t_u^+u)=(t_u^+)^2\|u\|^2+(t_u^+)^p\lambda \int_{\varOmega }|x|^{-s_{1}}|u|^pdx-(t_u^+)^q\int_{\varOmega }|x|^{-s_{2}}|u|^qdx=(t_u^+)^2g(t_u^+)=0.$$
Then
\begin{eqnarray*}
	\psi(t_u^+u)&=&(t_u^+)^2\|u\|^2+(p-1)(t_u^+)^p \lambda \int_{\varOmega }|x|^{-s_{1}}|u|^pdx-(q-1)(t_u^+)^q\int_{\varOmega }|x|^{-s_{2}}|u|^qdx-\varphi(t_u^+u)\\
	&=&(p-2)(t_u^+)^p \lambda \int_{\varOmega }|x|^{-s_{1}}|u|^pdx-(q-2)(t_u^+)^q\int_{\varOmega }|x|^{-s_{2}}|u|^qdx\\
	&=&(t_u^+)^p \frac{2pq}{p-q}\left( m^+-\frac{m}{2q}\right) -(t_u^+)^q\frac{2pq}{p-q} \left( m^+-\frac{m}{2p}\right) \\
\end{eqnarray*}
by (\ref{10}) and (\ref{11}).
Consider the function
\begin{eqnarray*}h(t)&\stackrel{\triangle}{=}&\frac{2pq}{p-q}\left(m^+-\frac{m}{2q}\right)t^p-\frac{2pq}{p-q}\left(m^+-\frac{m}{2p}\right)t^q,
\end{eqnarray*}
we have
\begin{eqnarray*}h'(t)&=&\frac{2p^2q}{p-q}\left(m^+-\frac{m}{2q}\right)t^{p-1}-\frac{2pq^2}{p-q}\left(m^+-\frac{m}{2p}\right)t^{q-1}.
\end{eqnarray*}
Hence, $h$ is strictly decreasing on $(0, t_0)$ and strictly increasing on $(t_0, +\infty)$, where $t_0=\left(\frac{q(m^+-\frac{m}{2p})}{p(m^+-\frac{m}{2q})}\right)^\frac{1}{p-q}$,
which implies that $\psi(t_u^+u)=h(t_u^+)<0$ by $h(0)=0$ and
\begin{eqnarray*}h(1)&=&\frac{2pq}{p-q}\left(m^+-\frac{m}{2q}\right)-\frac{2pq}{p-q}\left(m^+-\frac{m}{2p}\right)=-m \leq 0.
\end{eqnarray*}
Then $t_u^+ u\in M^+$.

Furthermore, we obtain
\begin{eqnarray*}
I(t_u^+u)
&=&\frac{1}{2}(t_u^+)^2\|u\|^2+\frac{1}{p}\lambda (t_u^+)^p\int_{\varOmega }|x|^{-s_{1}}|u|^pdx-\frac{1}{q}(t_u^+)^q\int_{\varOmega }|x|^{-s_{2}}|u|^qdx-\frac{1}{2}\varphi(t_u^+u)\\
&=&-\frac{p-2}{2p}\lambda (t_u^+)^p\int_{\varOmega }|x|^{-s_{1}}|u|^pdx+\frac{q-2}{2q}(t_u^+)^q\int_{\varOmega }|x|^{-s_{2}}|u|^qdx\\
&=&-\frac{q}{p-q}(t_u^+)^p\left(m^+-\frac{m}{2q}\right)+\frac{p}{p-q}(t_u^+)^q\left(m^+-\frac{m}{2p}\right)
\end{eqnarray*}
by (\ref{10}) and (\ref{11}).

Consider the function
\begin{eqnarray*}h(t)&\stackrel{\triangle}{=}&-\frac{q}{p-q}t^p\left(m^+-\frac{m}{2q}\right)+\frac{p}{p-q}t^q\left(m^+-\frac{m}{2p}\right),
\end{eqnarray*}
we have
\begin{eqnarray*}h'(t)&=&-\frac{pq}{p-q}t^{p-1}\left(m^+-\frac{m}{2q}\right)+\frac{pq}{p-q}t^{q-1}\left(m^+-\frac{m}{2p}\right).
\end{eqnarray*}
Hence, $h$ is strictly increasing on $(0, t_0)$, where $t_0=\left(\frac{m^+-\frac{m}{2p}}{m^+-\frac{m}{2q}}\right)^\frac{1}{p-q}\geq1$,
which implies that
$$
m^+\leq I(t_u^+u)=h(t_u^+)<h(1)=m^+.
$$
That is $\lim\limits_{n\to\infty}\|u_n\|=\|u\|$. Then we have that $u_n \to u$ in $\HO$. $I(u)=m^+$ and $u \in M^+ \cup M^0$. If $u \in M^0$, $M^0$ is nonempty. From Proposition \ref{p5} , $m^+<m^0$ , but $m^0 \leq I(u)=m^+$, so $u \in M^+$. Hence
$\nabla(I|_{M^+})(u)=0$. Thus, there exists $\lambda \in \RR$ such that $$\nabla I(u)=\nabla(I|_{M^+})(u)+\lambda\nabla \varphi(u)=\lambda\nabla \varphi(u),$$ which implies that  $$0=\varphi(u)=(\nabla I(u), u)=\lambda(\nabla\varphi(u),u)=\lambda \psi(u).$$ Then  $\lambda=0$ by $\psi(u)<0$, so $\nabla I(u)=0$, $I'(u)=0$.
 Then $u$ is a solution of \cref{eq1}. Note that $u\geq 0$ and $u\neq 0$. By the strong maximum principle, $u$ is a positive solution of \cref{eq1}. Moreover, it follows from Lemma 2.6 in \cite{Tang2025} that $u\in C(\overline{\Omega})$, which completes the proof. $\hfill\Box$

{\bf Proof of Theorem \ref{thm2}}\ \ \ \
We assume that $u$ is a solution of \cref{eq1}, then the $Pohozave$ identity is
\begin{eqnarray*}
	\frac{1}{2}\int_{\partial\varOmega}  x\cdot\nu \left( \dfrac{\partial u}{\partial \nu}\right)^2dS+\frac{N-2}{2}\|u\|^2+\frac{\lambda (N-s_1)}{p}\int_{\varOmega }|x|^{-s_{1}}|u|^pdx-\frac{N-s_2}{q}\int_{\varOmega }|x|^{-s_{2}}|u|^qdx=0,
\end{eqnarray*}
where $\nu$ denotes the outward normal to $\partial\varOmega$. The $Nehari$ identity is
\begin{eqnarray*}
\|u\|^2+\lambda\int_{\varOmega }|x|^{-s_{1}}|u|^pdx-\int_{\varOmega }|x|^{-s_{2}}|u|^qdx=0.
\end{eqnarray*}
When $q=2^*(s_2)$, we have
\begin{eqnarray*}
	\frac{1}{2}\int_{\partial\varOmega}  x\cdot\nu \left( \dfrac{\partial u}{\partial \nu}\right)^2dS+\frac{N-2}{2}\|u\|^2+\frac{\lambda (N-s_1)}{p}\int_{\varOmega }|x|^{-s_{1}}|u|^pdx-\frac{N-2}{2}\int_{\varOmega }|x|^{-s_{2}}|u|^{2^*(s_2)}dx=0
\end{eqnarray*}
and
\begin{eqnarray*}
\|u\|^2+\lambda\int_{\varOmega }|x|^{-s_{1}}|u|^pdx-\int_{\varOmega }|x|^{-s_{2}}|u|^{2^*(s_2)}dx=0,
\end{eqnarray*}
which implies that
\begin{eqnarray*}
	\frac{1}{2}\int_{\partial\varOmega}  x\cdot\nu \left( \dfrac{\partial u}{\partial \nu}\right)^2dS+\left( \frac{N-s_1}{p}-\frac{N-2}{2}\right) \lambda\int_{\varOmega }|x|^{-s_{1}}|u|^pdx=0.
\end{eqnarray*}
We note that $\ \frac{N-s_1}{p}-\frac{N-2}{2}>0$ from $p<2^*(s_1)$ and $ x\cdot\nu \geq 0$ for $x \in \partial\varOmega$, because that $\varOmega$ is a starshaped domain about the origin. Then $u\equiv 0$, \cref{eq1} has no nonzero solution.$\hfill\Box$

{\bf Proof of Theorem \ref{thm3}}\ \ \ \ Let $m^+_q=m^+$, $I_q(u)=I(u)$ and $\psi_q(u)=\psi(u)$ with $2<q\leq 2^*(s_2)$. Because that $u_{q_j}$ is the positive solution of \cref{eq1} by \cref{01}, we have
\begin{align*}
	I_{q_j}'(u_{q_j})=0,\label{I`=0}\\
	I_{q_j}(u_{q_j})=m_{q_j}^+,\\
	\psi_{q_j}(u_{q_j})\leq0.
\end{align*}
We know $\limsup\limits_{j \to \infty} m_{q_j}^+\leq m_{2^*(s_2)}^+$ from Proposition \ref{p6}. By  (\ref{eq2.2}), for $j$ large enough we have
$$m_{2^*(s_2)}^+ +1 \geq m_{q_j}^+=I_{q_j}(u_{q_j})>\frac{(p-2)(q_j-2)}{2pq_j}\|u_{q_j}\|^2\geq \frac{(p-2)(q_1-2)}{2p{2^*(s_2)}}\|u_{q_j}\|^2,$$
which implies that $\left\lbrace u_{q_j}\right\rbrace $ is bounded in $\HO$. Moreover, $\left\lbrace u_{q_j}\right\rbrace $ is a PS sequence of $I_{2^*(s_2)}(u)$. Going if necessary to a subsequence, we can assume that
\begin{eqnarray*}
	u_{q_j}\rightharpoonup u \ \ \ \ \  \ \ \ \ \ \  \ \ \  \ \text{in}\ H^1_0{(\varOmega)},
\end{eqnarray*}
Then $I'_{2^*(s_2)}(u)=0$. From \cref{thm2}, we know $u=0$ and $u_{q_j} \rightharpoonup 0$.

We firstly note
\begin{eqnarray*}\label{mqj0}
	M_{q_j}\stackrel{\bigtriangleup}{=}u_{q_{j}} ( x_{q_{j}} )\stackrel{\bigtriangleup}{=}\operatorname* {m a x}\limits_{x \in \varOmega} u_{q_{j}}(x) \to\infty
\end{eqnarray*}
as $j \to\infty$. Otherwise, there exists $M>1$ such that $u_{q_j}(x)\leq M_{q_j}\leq M$ for $x \in \varOmega$, and
\begin{eqnarray*}
	\int_{\varOmega} {| x |^{-s_2}}{u_{q_j}^{2^{*} ( s_2 )-\varepsilon_j}}dx&	=&\int_{\varOmega} {| x |^{-s_2}}{u_{q_j}^{2^{*} ( s_2 )-1+1-\varepsilon_j}} dx\\
	&\leq&M^{1-\varepsilon_j}\int_{\varOmega} {| x |^{-s_2}}{u_{q_j}^{2^{*} ( s_2 )-1}} dx\\
	&\leq&M\int_{\varOmega}{| x |^{-s_2}}{u_{q_j}^{2^{*} ( s_2 )-1}} dx\\
	&\rightarrow&0
\end{eqnarray*}
as $j\to \infty$, due to $u_{q_j} \rightharpoonup 0$ and the compactness of $\HO \hookrightarrow L^{2^*(s_2)-1}(|x|^{-s_2};\varOmega)$. Similarly, one has
\begin{eqnarray*}
	\int_{\varOmega} {| x |^{-s_1}}{u_{q_j}^{p}}dx \rightarrow 0~~as~~j\to \infty.
\end{eqnarray*}
From the $Nehari$ identity
\begin{eqnarray*}
	\|u_{q_j}\|^2+\lambda\int_{\varOmega }|x|^{-s_{1}}|u_{q_j}|^pdx-\int_{\varOmega }|x|^{-s_{2}}|u_{q_j}|^{2^{*} ( s_2 )-\varepsilon_j}dx=0
\end{eqnarray*}
we obtain $||u_{q_j}|| \to 0$ as $j\to \infty$. Then $u_{q_{j}} \to0$ strongly in $H_{0}^{1} ( \varOmega)$. It follows from (\ref{2.0}) and Hardy-Sobolev inequality that
\begin{eqnarray*}
\|u_{q_{j}}\|^2&\leq&\|u_{q_{j}}\|^2+\lambda\int_{\varOmega }|x|^{-s_{1}}|u_{q_{j}}|^pdx\\
&=&\int_{\varOmega }|x|^{-s_{2}}|u_{q_{j}}|^{2^*(s_2)-\varepsilon_{j}}dx\\
&=&\int_{\varOmega }|x|^{-s_{2}}|u_{q_{j}}|^{q_{j}}dx\\
&\leq&M^{q_j-q_1}\int_{\varOmega}{| x |^{-s_2}}{u_{q_j}^{q_1}} dx\\
&\leq&M^{\varepsilon_1}\int_{\varOmega}{| x |^{-s_2}}{u_{q_j}^{q_1}} dx\\
&\leq&C\|u_{q_{j}}\|^{q_1}
\end{eqnarray*}
for some $C>0$. Then we have $||u_{q_j}||\geq C^{-\frac{1}{q_1-2}}$, which is a contradiction.

Set
\begin{eqnarray*}\label{vqj}
v_{q_{j}} ( y ) :=\frac{u_{q_{j}} ( x_{q_{j}}+\kappa_{j} y )} {M_{q_j}} \quad\mathrm{for}~y\in\varOmega_j,
\end{eqnarray*}
where
\begin{eqnarray*}\label{2.7.1}
	\kappa_{j} :=M_{q_j}^{-\frac{2}{N-2}},
\end{eqnarray*}
and
\begin{eqnarray*}\label{oj}
\varOmega_{j} :=\big\{y \in\mathbb{R}^{N} \, \big| \, x_{q_{j}}+\kappa_{j} y \in\varOmega\big\}.
\end{eqnarray*}

Suppose that $x_{q_j}\to x_0 \in \overline{\Omega}$. Assume that $x_0 \in {\Omega}$. Let $y\in B_r(0)$, then $\kappa_{j}y+x_{q_j} \to x_0$ as $j\to \infty$, which implies that $\kappa_{j}y+x_{q_j}\in \overline{B}_\delta(x_0)\subset \Omega$. Then we get $y\in \Omega_j$, that is, $B_r(0) \subset \Omega_j$ for $j$ large enough.

Moreover, we have
\begin{eqnarray*}
	\int_{\varOmega }|\nabla u_{q_j}|^2dx&=&{M_{q_j}^2}{\kappa_{j}^{N-2}}\int_{\varOmega_j }|\nabla v_{q_{j}}|^2dy=\int_{\varOmega_j }|\nabla v_{q_{j}}|^2dy,
\end{eqnarray*}
Then $v_{q_j}\in H^1_0(\varOmega_{j})\subset D^{1,2}(\mathbb{R}^{N})$, $\{v_{q_j}\}$ is a bounded sequence of $D^{1,2}(\mathbb{R}^{N})$. Going if necessary to a subsequence, we can assume that
\begin{eqnarray*}
	&&v_{q_j}\rightharpoonup v \ \ \ \ \  \ \ \ \ \ \  \ \ \  \ \text{in}\ D^{1,2}{(\mathbb{R}^{N})},\\
	&&v_{q_j}\rightarrow v \ \ \ \ \ \ \ \ \  \ \ \ \ \ \ \text{in}\ L^{q}_{loc}(\mathbb{R}^{N})\ (q\in [2,2^*)),\\
	&&v_{q_j}\rightarrow v(x) \ \ \ \ \  \ \ a.e. \ \text{in}\  \mathbb{R}^{N},
\end{eqnarray*}
as $j\rightarrow\infty$. Moreover, $v\geq 0$ and  $\lim\limits_{j\to\infty}\|v_{q_j}\|\geq\|v\|$.

Next we prove that $v_{q_j} \rightharpoonup 0$ as $j\to \infty$. Choose a $w\in C_0^1(\mathbb{R}^{N})$, we only need to prove that $\int_{\mathbb{R}^{N}}\nabla v_{q_j} \cdot \nabla w\to 0$. Let

$$
w_{j}(x) :={\kappa_{j}}^{-\frac{N-2}{2}}w\left(\frac{x-x_{q_{j}}}{\kappa_{j}}\right).
$$
Then
\begin{eqnarray*}
	\int_{\mathbb{R}^{N}}|\nabla w_{j}|^2dx&={\kappa_{j}^{-N}}\int_{\mathbb{R}^{N}}|\nabla w\left(\frac{x-x_{q_{j}}}{\kappa_{j}}\right)|^2dx=\int_{\mathbb{R}^{N}}|\nabla w(y)|^2dy
\end{eqnarray*}
and
\begin{eqnarray*}
	\int_{\mathbb{R}^{N}}\nabla v_{q_j}(y) \cdot \nabla w(y) dy &=&{M_{q_j}^{-1}}{\kappa_{j}}\int_{\mathbb{R}^{N}}\nabla u_{q_j}(x_{q_{j}}+\kappa_{j} y) \cdot \nabla w(y) dy\\
	&=&{M_{q_j}^{-1}}{\kappa_{j}^{1-N}}\int_{\mathbb{R}^{N}}\nabla u_{q_j}(x) \cdot \nabla w\left(\frac{x-x_{q_{j}}}{\kappa_{j}}\right)dx\\
	&=&\int_{\mathbb{R}^{N}}\nabla u_{q_j}(x) \cdot \nabla w_j(x)dx.
\end{eqnarray*}
Because that $I_{q_j}'(u_{q_j})=0$, one has
\begin{eqnarray*}
	0&=&\langle I_{q_j}'(u_{q_j}),\ w_j \rangle \\
	&=&\int_{\mathbb{R}^{N}}\nabla u_{q_j} \cdot \nabla w_j dx+\lambda\int_{\mathbb{R}^{N} }|x|^{-s_{1}}|u_{q_j}|^{p-2}u_{q_j}w_j dx- \int_{\mathbb{R}^{N}}|x|^{-s_{2}}|u_{q_j}|^{q_j-2} u_{q_j}{w_j}dx\\
	&=&\int_{\mathbb{R}^{N}}\nabla v_{q} \cdot \nabla w dx+\lambda\int_{\mathbb{R}^{N} }|x|^{-s_{1}}|u_{q_j}|^{p-2}u_{q_j}w_j dx- \int_{\mathbb{R}^{N}}|x|^{-s_{2}}|u_{q_j}|^{q_j-2} u_{q_j}{w_j}dx.
\end{eqnarray*}
We note that
\begin{eqnarray*}
	\left| \int_{\mathbb{R}^{N}}|x|^{-s_{2}}|u_{q_j}|^{q_j-2} u_{q_j}{w_j}dx\right| &=&
	{M_{q_j}^{q_j-1}}{\kappa_{j}^{\frac{N}{2}+1}}\left| \int_{\mathbb{R}^{N}}|x_{q_j}+\kappa_{j}y|^{-s_{2}}|v_{q_j}|^{q_j-2} v_{q_j}{w(y)}dy\right| \\
	&\leq&
	\frac{\kappa_{j}^{s_2}}{|x_{q_{j}}|^{s_2}}{M_{q_j}^{q_j-1}}{\kappa_{j}^{\frac{N}{2}+1-s_2}}\int_{\mathbb{R}^{N}}\left|  \frac{x_{q_{j}}}{|x_{q_{j}}|}+\frac{\kappa_{j} y}{|x_{q_{j}}|}\right|^{-s_{2}}|v_{q_j}|^{q_j-2} v_{q_j}{|w(y)|}dy\\
	&=&
	\frac{\kappa_{j}^{s_2}}{|x_{q_{j}}|^{s_2}}{M_{q_j}^{q_j-{2^*(s_2)}}}\int_{\mathbb{R}^{N}}\left|  \frac{x_{q_{j}}}{|x_{q_{j}}|}+\frac{\kappa_{j} y}{|x_{q_{j}}|}\right|^{-s_{2}}|v_{q_j}|^{q_j-2} v_{q_j}{|w(y)|}dy\\
	&\leq&		\frac{\kappa_{j}^{s_2}}{|x_{q_{j}}|^{s_2}}2^{s_2}\int_{\mathbb{R}^{N}}|v_{q_j}|^{q_j-2} v_{q_j}{|w(y)|}dy\\
		&\leq&		\frac{\kappa_{j}^{s_2}}{|x_{q_{j}}|^{s_2}}2^{s_2}\int_{\mathbb{R}^{N}}{|w(y)|}dy\\
	&\leq&		\frac{\kappa_{j}^{s_2}}{|x_{q_{j}}|^{s_2}}2^{s_2}||w||_{\infty}|\mathrm{supp}\ w|\\
	&\to&0
\end{eqnarray*}
as $j\to\infty$, where $|\mathrm{supp}\ w|$ is the measure of $\mathrm{supp}\ w$. Similarly, we can obtain that  $\lambda\int_{\mathbb{R}^{N} }|x|^{-s_{1}}|u_{q_j}|^{p-2}u_{q_j}w_j dx\to 0$ as $j\to\infty$. Then we prove that $v_{q_j} \rightharpoonup 0$ as $j\to \infty$.

Notice that $v_{q_j}(0)=1$, $0\leq v_{q_j}(x) \leq 1$ for $x \in \Omega_j$, and
$$
 {-\Delta v_{q_{j}}+\lambda\kappa_{j}^{2} M_{q_j}^{p-2}   \left( \frac{\kappa_{j}} {| x_{q_{j}} |} \right)^{s_1} \frac{v_{q_{j}}^{p-1}} {\left|  \frac{x_{q_{j}}} {| x_{q_{j}} |}+\frac{\kappa_{j}} {| x_{q_{j}} |} y \right| ^{s_1}}-\left( \frac{\kappa_{j}} {| x_{q_{j}} |} \right)^{s_2} \frac{v_{q_{j}}^{q_{j}-1}} {\left|  \frac{x_{q_{j}}} {| x_{q_{j}} |}+\frac{\kappa_{j}} {| x_{q_{j}} |} y \right| ^{s_2}}=0 \quad\mathrm{i n} \ \ B_r(0) (\subset \varOmega_{j})}.
$$
Hence, for every $q > N$ we obtain
$$ \left( \frac{\kappa_{j}} {| x_{q_{j}} |} \right)^{s_2} \frac{v_{q_{j}}^{q_{j}-1}} {\left|  \frac{x_{q_{j}}} {| x_{q_{j}} |}+\frac{\kappa_{j}} {| x_{q_{j}} |} y \right| ^{s_2}}\in L^q(B_r(0))$$
from
\begin{eqnarray*}
	\int_{B_r(0)}\left( \left( \frac{\kappa_{j}} {| x_{q_{j}} |} \right)^{s_2} \frac{v_{q_{j}}^{q_{j}-1}} {\left|  \frac{x_{q_{j}}} {| x_{q_{j}} |}+\frac{\kappa_{j}} {| x_{q_{j}} |} y \right| ^{s_2}}\right) ^qdy	&\leq& 2^{s_2}|B_r(0)|.
\end{eqnarray*}
Similarly, one has
$$ \lambda\kappa_{j}^{2} M_{q_j}^{p-2}   \left( \frac{\kappa_{j}} {| x_{q_{j}} |} \right)^{s_1} \frac{v_{q_{j}}^{p-1}} {\left|  \frac{x_{q_{j}}} {| x_{q_{j}} |}+\frac{\kappa_{j}} {| x_{q_{j}} |} y \right| ^{s_1}}\in L^q(B_r(0)).$$
Thus, elliptic theory implies $||v_{q_j}||_{W^{2,q}(B_r(0))}$ is bounded. Then we obtain $||v_{q_j}||_{C^{1}(\overline{B_r(0)})}$ is bounded by Sobolev's imbedding theorem. Then $v_{q_j}$
 is equicontinuous on $\overline{B_r(0)}$, hence by Arzela-Ascoli theorem, there exists a subsequence converging to some $V$ uniformly on $\overline{B_r(0)}$, which contradicts $v_{q_j} \rightharpoonup 0$ as $j\to \infty$ due to $V(0)=1$. Then $x_0 \in \partial\Omega.$$\hfill\Box$

\bibliographystyle{unsrt}
\bibliography{ref}

\end{document}